\documentclass[
submission
 , nomarks
]{dmtcs-episciences}

\usepackage[utf8]{inputenc}
\usepackage{subfigure}

\usepackage{enumerate}

\usepackage{amsmath,amssymb}

\newtheorem{defi}{Definition}[section]

\newtheorem{rmk}[defi]{Remark}

\newtheorem{prop}[defi]{Proposition}
\newtheorem{theo}[defi]{Theorem}
\newtheorem{cor}[defi]{Corollary}

\newcommand{\Card}{\text{Card}}

\newcommand{\D}{\mathcal{D}}
\newcommand{\X}{\mathcal{X}}

\usepackage[round]{natbib}

\author{R\'emi Molinier\affiliationmark{1}}

\title[Comparing $n$-dice]{Comparing $n$-dice fixing the sum of the faces}
\affiliation{Univ. Grenoble Alpes, CNRS, IF, 38000 Grenoble, France}
 
\begin{document}
%\publicationdetails{VOL}{2015}{ISS}{NUM}{SUBM}
\maketitle
\begin{abstract}
These notes describe some results on dice comparisons when changing the numbers on the faces while the sum of all the faces stay the same. 
\end{abstract}

\section{Introduction}
\label{sec:in}

These notes describe some results on dice comparisons when changing the numbers on the faces. This project started as an attempt to create an activity for children about probability and dice. It was motivated by the existence of non-transitive dice (see for example \citep{BGH,AD,S}). Apart from the fun part of such dice, Their existence and more generally dice games \citep{T2} can illustrate Arrow's impossibility Theorem \citep{A}. \citet {T3} also gave some nice asymptotic results on the proportion of ties in some families of dice.

The idea was to elaborate some kind of tournament where each player creates its own die by choosing the number on each face. They then play against each other rolling their dice. With no other rules, it is obvious that the best way to win would be to be the player choosing the highest numbers on the faces of its die. Therefore, some conditions have to be set on the faces to give some interesting features to the game. The author idea, and the point of view in these notes, is to only allow positive or non negative numbers on the faces and to fix the value of the sum of all the faces (or equivalently the mean of all the faces). One natural question arise then, is there a "best" die.  Here we are comparing two dice $A$ and $B$ by looking at the probabilities that $A$ rolls higher than $B$ and that $B$ rolls higher than $A$.  The existence of a "best" die is not always guaranteed  as one can check on small dice with small fixed sum value. However, when we fix the sum  to be $1+2+3+4+5+6=21$, i.e. the sum of the face of the standard die $D_{\text{st}}=(1,2,3,4,5,6)$, then you cannot beat the standard die $D_{\text{st}}$ if you only allows positive numbers. You can even see that $D_{\text{st}}$ is always better except when all your faces are in between 1 and 6 included; In that case, you will actually be as good as $D_{\text{st}}$. Hence you need a zero to beat $D_{\text{st}}$ and we can actually see how to construct your die to be sure to beat it. 

This paper presents some general consideration along these lines.  We work here with $n$-dice, i.e. dice with $n$ faces and we prove some general result about comparing $n$-dice when fixing the sum of the faces. For example we prove the the $n$-die $(0,1,2,3,\dots,n-1)$ cannot be beaten by another $n$-die where the sum of the faces is $n(n-1)/2$. As noticed before, there is not a "best" die in general. However, there is always a worst die (still fixing the sum of the faces) which is one with 0 on each faces except one. Some aspects of these results were already, and independently, known and studied by\citet{T1}. In \citep{T1}, the author study families of dice, fixing the maximum value, minimum value and the sum of the faces, where there are a die which ties with all the over dice of the family. In particular, The author give in \citep[Theorem 2]{T1} a nice characterization of such a die. 

In the last section, we compare $n$-dice, with face sum equal to $n(n+1)/2$ to the standard $n$-die. We give a characterization of such $n$-dice which beat, loose or tie with the standard $n$-die, in terms of the numbers on the faces and give then ones which has the give the higher probability to beat the standard $n$-die.
\medskip

\noindent \textbf{Organization.} The paper is organized as follow. Section \ref{sec:def} fixes the notations and definitions. In Section \ref{sec:unbeatable} we explain why the die $(0,1,2,\dots,n-1)$ is unbeatable in its class of $n$-dice and give some corollaries. We also look at the worst die. Finally, in Section \ref{sec:usual} we compare $n$-dice with the standard $n$-die $D_{n,\text{st}}=(1,2,\dots,n)$.   
\medskip

\noindent \textbf{ Acknowledgments.} I would like to thank Dave Auckly with whom I started the idea of creating an activity around these ideas. This work wouldn't be achieved without enlightening discussions with Sylvain Gravier and Florian Galliot and. I am also grateful to Emilie Devijver for her support, some useful discussion on the subject, and some help on some computational aspects. Finally, I would like to thank Lorenzo Traldi for, after the first version of these notes were put on Arxiv, reaching out and enlightening me with references on the subject.  

% You may scarsely use \clearpage to advance to a new page if this
% improves the readability of the document structure
%\clearpage

\section{Definitions and notations}
\label{sec:def}

Let $n$ be an integer greater than 2.

\begin{defi} 
 A \emph{$n$-die} is an increasing sequence of non-negative integers $D=(f_1,f_2,\dots,f_n)$ and, for $i\in\{1,2,\dots,n\}$ $f_i$ is the \emph{$i$th face} of $D$. We denote by $\D_n$ the set of all $n$-dice.

If $\sigma\in\natural$, a \emph{$(\sigma,n)$-die} is a $n$-die $D=(f_1,f_2,\dots,f_n)$ such that 
\[
 \sum_{i=1}^n f_i=\sigma.
\]
We denote by $\D_n(\sigma)$ the set of all $(\sigma,n)$-dice. 
\end{defi}

\begin{defi}\label{defi:relation}
 Let $D=(f_1,f_2,\dots,f_n)$ and $D'=(f_1',f_2',\dots,f_n')$ be two dice.
 
 We denote by $\gamma(D,D')$, respectively $\eta(D,D')$ the number of time $D$ rolls higher than $D'$, resp. $D$ is equal to $D'$, when looking at all possible issues of rolling the two dice at the same time. In other words,
\begin{align*}
 \gamma(D,D')&=\sum_{i=1}^{n}\Card\left\{j\,\mid\,f_i>f_j'\right\}\\
 \eta(D,D')&=\sum_{i=1}^{n}\Card\left\{j\,\mid\,f_i=f_j'\right\}
\end{align*}
 
We then denote by $\Delta(D,D')=\gamma(D,D')-\gamma(D',D)$ the differential between the two dice. 

Finally, we set the following notations
\begin{align*}
 D \succ D' & \text{ if }\Delta(D,D')>0\\
 D \prec D' & \text{ if }\Delta(D,D')<0\\
 D\sim D' &\text{ if }\Delta(D,D')=0\\
 D \succsim D' &\text{ if }\Delta(D,D')\geq0\\
 D \precsim D' &\text{ if }\Delta(D,D')\leq0\\
\end{align*}

\end{defi}

Be aware that, in general, $\prec$ is \underline{not} a partial order on $\D_n$ or even $\D_n(\sigma)$ when $\sigma>2$. The next proposition gives some easy properties.

\begin{prop}\label{prop:basic}
 Let $D$ and $D'$ be two $n$-dice.
 \begin{enumerate}[(a)]
  \item $0\leq\gamma(D,D')\leq n^2$,
  $0\leq\eta(D,D')\leq n^2$, and
  $-n^2\leq\Delta(D,D')\leq n^2$.
  \item $\eta(D,D')=\eta(D',D)$ and $\Delta(D,D')=-\Delta(D',D)$.
  \item $\gamma(D,D')+\gamma(D',D)+\eta(D,D')=n^2$.
 \end{enumerate}
\end{prop}

\begin{proof}
 These properties are easily derived from Definition \ref{defi:relation} and the fact that there is $n^2$ comparisons between the faces of $D$ and the faces of $D'$.
\end{proof}

\section{An unbeatable die and the worst die} %in ${\D_n\left(\frac{n(n-1)}{2}\right)}$}
\label{sec:unbeatable}

Let $D_0=(0,1,2,3,\dots,n-1)$. We have $D_0\in{\D_n\left(\frac{n(n-1)}{2}\right)}$ and we are going to see that $D_0$ is unbeatable in $ \D_n\left(\frac{n(n-1)}{2}\right)$.

\begin{prop}\label{prop:D0}
 Let $D=(f_1,f_2,\dots,f_{n-1})\in{\D_n\left(\frac{n(n-1)}{2}\right)}$.
 
 Then $D\precsim D_0$ and $D\sim D_0$ if and only if, for all $i\in\{1,2,\dots,n\}$, $f_i\leq n-1$.
 
\end{prop}

\begin{proof}
 In this context, we have 
 \[
 \gamma(D,D_0)= \sum_{i=1}^n \min(f_i,n).
 %\qquad\text{and}\qquad \gamma(D_0,D)=\sum_{i=1}^n \max(n-1-f_i,0)
 \]
 In particular, 
 \[
  \gamma(D,D_0)\leq \sum_{i=1}^n f_i=\frac{n(n-1)}{2}.
 \]
 Moreover, 
 \[
\eta(D_0,D)=\eta(D,D_0)=\sum_{i=0}^{n-1}\Card\{j\,\mid\, f_j=i\}=\Card\{j\,\mid\, f_j\leq n-1\}\leq n
 \]
 Hence, by Proposition \ref{prop:basic}, 
 \[
  \gamma(D_0,D)=n^2-\gamma(D,D_0)-\eta(D,D_0)\geq n^2- \frac{n(n-1)}{2}-n=\frac{n(n-1)}{2}.
 \]
 In particular, $\gamma(D_0,D)\geq\gamma(D,D_0)$.
 
 Moreover, we have $\gamma(D_0,D)=\gamma(D,D_0)$ if and only if 
 \[
  \gamma(D,D_0)=\frac{n(n-1)}{2} \qquad\text{and}\qquad \eta(D,D_0)=n,
 \]
 i.e., for all $i\in\{1,2,\dots,n\}$, $f_i\leq n-1$.
\end{proof}

\begin{cor}\label{cor:main}
 Let $p\geq 0$, $q=p+n-1$ and $\sigma=\sum_{i=p}^{q} i=\frac{n(p+q)}{2}$. 
 Set $D_{p,q}=(p,p+1,\dots,q)\in\D_n(\sigma)$ and let $D=(f_1,f_2,\dots,f_n)\in \D_n(\sigma)$.
 
 If, for all $i\in\{1,2,\dots,n\}$, $f_i\geq p$, then $D\precsim D_{p,q}$.
 Moreover, $D\sim D_{p,q}$ if and only if, for all $i\in\{1,2,\dots,n\}$, $p\leq f_i\leq q$. 
\end{cor}

\begin{proof}
 Let 
 \begin{align*}
  \widetilde{D}=(f_1-p,f_2-p,\dots,f_n-p).
 \end{align*}
 Then we have that $\Delta(D,D_{p,q})=\Delta(\widetilde{D},D_{0,n-1})$ and the result follows from Proposition \ref{prop:D0}.
\end{proof}

The fact that, for any $n$-die $ D\in \D_n\left(\frac{n(n-1)}{2}\right)$ with faces in $\{0,\dots,n-1\}$ we have $D\sim D_0$, was already known by Traldi and is a consequence of \citep[Corollary 5]{T1}.
Florian Galliot also suggested an other nice argument. First you need to notice the following fact.  Let $\X_n$ be the set of all $n$-die $D=(f_1,f_2,\dots,f_n)\in\D_n\left(\frac{n(n-1)}{2}\right)$ such that for all $i\in\{1,2,\dots,n\}$, $0\leq f_i\leq n-1$. 

\begin{prop}\label{prop:swap}
Let $D=(f_1,f_2,\dots,f_n)\in \X_n$.
Let $i_0,j_0\in\{1,\dots,n\}$ and set \[\widetilde{D}=(f_1,\dots,f_{i_0}+1,\dots,f_{j_0}-1,\dots,f_n).\]

If $f_{i_0}\leq n-2$ and $f_{j_0}\geq 1$ (or, equivalently, $\widetilde{D}\in \X_n$), then
\[\Delta(D,D_0)=\Delta(\widetilde{D},D_0).\]
\end{prop}

Notice that the conditions,  $f_{i_0}\leq n-2$ and $f_{j_0}\geq 1$.

\begin{proof}
 We have 
 \begin{align*}
\gamma(D,D_0)&= \sum_{i=1}^n \min(f_i,n)\\
&=f_{i_0}+f_{j_0}+\sum_{i\in\{1,2,\dots,n\}\smallsetminus \{i_0,j_0\}} \min(f_i,n)\\
&= (f_{i_0}+1)+(f_{j_0}-1)+\sum_{i\in\{1,2,\dots,n\}\smallsetminus \{i_0,j_0\}} \min(f_i,n)\\
&=\gamma(\widetilde{D},D_0),  
 \end{align*}
 and $\eta(D,D_0)=n=\eta(\widetilde{D},D_0)$.
 Thus, by Proposition \ref{prop:basic}, we also have $\gamma(D_0,D)=\gamma(D_0,\widetilde{D})$. Therefore $\Delta(D,D_0)=\Delta(\widetilde{D},D_0)$. 
\end{proof}

\begin{cor}\label{cor:1}
Let $D=(f_1,f_2,\dots,f_n)\in  \D_n\left(\frac{n(n-1)}{2}\right)$ be  such that for all $i$, $0\leq f_i\leq n-1$.

Then $D\sim D_0$.
\end{cor}

\begin{proof}
One can notice that, for any $D\in \X_n$ there exists a sequence of dice \[D_0=X_0,X_1,X_2,\dots,X_r=D\] such that for all $i\in\{1,2,\dots,r-1\}$, $X_i\in \X_n$ and that $X_{i+1}$ is obtained from $X_i$ by adding 1 on one face and subtracting 1 to another one (for example, one can first change the highest face till it become $n$, then the second highest till it becomes $n-1$, and so on). 
%To do so, we can prove, by induction on $n$, that any die in $\X_n$ can be linked to $(0,1,\dots,n-1)$. For $n=1$ it is trivial and if we assume that it is true for $n-1$, then we can first ... bla bla bla 
Then, by the Proposition \ref{prop:swap}, $\Delta(D,D_0)=\Delta(D_0,D)=0$.   
\end{proof}

%Actually, Corollary \ref{cor:1} was already known by Traldi and is also a corollary of \citep[Theorem 2]{T1} and it give a third different proof.\\

We finish this section by giving the worst die. Notice that the existence of such a die doesn't depends on the fixed  value for the sum of the faces.

\begin{prop}
 Let $\sigma\geq 2$. 
 Set $D_w=(0,0,\dots,0,\sigma)\in\D_n(\sigma)$ and let $D\in\D_n(\sigma)\smallsetminus \{D_w\}$. 
 
 If $n\geq 3$, then $D\succ D_w$. 
\end{prop}

\begin{proof}
 set $D=(f_1,f_2,\dots,f_n)$.
 Since $D\in \D_n(\sigma)\smallsetminus \{D_w\}$, $0<f_n<\sigma$ and $0<f_{n-1}<\sigma$. Thus $\gamma(D,D_w)\geq 2(n-1)=2n-2$. Moreover, $\gamma(D_w,D)=n$. Thus $\Delta(D,D_w)\geq (2n-2)-n=n-2>0$.
\end{proof}

\section{Comparison with the standard die}\label{sec:usual}

In this section, we are interested in comparing $n$-dice with the \emph{standard $n$-die} $D_{n,\text{st}}$ (or $D_{\text{st}}$ when $n$ is understood) given by $D_{n,\text{st}}=(1,2,3,\dots,n).$. Notice that $D_{n,\text{st}}\in \D_n\left(\frac{n(n+1)}{2}\right)$. By Corollary \ref{cor:main}, we know that $D_{n,\text{st}}$ can not be beaten by a die $D\in\D_n\left(\frac{n(n+1)}{2}\right)$ with only positive faces but we can be more precise. We give here a characterization of the $n$-dice $D\in\D_n\left(\frac{n(n+1)}{2}\right)$ such that $D\prec D_{n,\text{st}}$ or $D\succ D_{\text{st}}$ and $D\sim D_{\text{st}}$ in terms of the faces of $D$.

\begin{prop}\label{prop:characterization}
 Let ${D\in\D_n\left(\frac{n(n+1)}{2}\right)}$ and let $k,l,r\in\naturals$ such that $D$ is given by
 \begin{equation}
 \label{eq:form}
 D=(\underbrace{0,0,\dots,0}_{k\text{ zeros}},\underbrace{f_1,f_2,\dots,f_l}_{\forall i,\; 1\leq f_i\leq n},\underbrace{g_1,g_2,\dots,g_r}_{\forall j,\; g_j\geq n+1}).
 \end{equation}

 Then, 
 \begin{align*}
 \frac{1}{2}\Delta(D,D_{\text{st}})= \left(r(n+1) - \sum_{i=1}^r g_i \right) +\frac{k-r}{2}.
 \end{align*}
\end{prop}

\begin{proof}
 We have 
 \begin{align*}
  &\gamma(D,D_{\text{st}})=rn + \sum_{i=1}^l (f_i-1)=rn-l+ \sum_{i=1}^l f_i,\quad \text{and}\\
  &\gamma(D_{\text{st}},D)=kn + \sum_{i=1}^l (n-f_i)=kn+ln - \sum_{i=1}^l f_i.
 \end{align*}
 Hence
 \begin{align*}
  \Delta(D,D_{\text{st}})&=\gamma(D,D_{\text{st}})-\gamma(D_{\text{st}},D)\\
  &=\left(rn-l+ \sum_{i=1}^l f_i\right)-\left(kn+ln - \sum_{i=1}^l f_i\right) &\\
  &=2\sum_{i=1}^l f_i -  n(k+l) + rn  -l  &\\
  &=2\sum_{i=1}^l f_i -  n(n-r)+ rn - (n-r-k)  & \text{(because }k+l+r=n)\\
  &=2\sum_{i=1}^l f_i -  (n^2+n)+ 2rn + r + k &\\
  &=2\sum_{i=1}^l f_i - n(n+1)+2r(n+1) -r+k&\\
  &=2\left(\sum_{i=1}^l f_i - \frac{n(n+1)}{2}+r(n+1)+\frac{k-r}{2}\right)&\\
   &=2\left(r(n+1) - \sum_{i=1}^r g_i +\frac{k-r}{2}\right)  & \left(\text{since }\sum_{i=1}^l f_i+\sum_{i=1}^r g_i = \frac{n(n+1)}{2}\right).
 \end{align*}
 and the result follows.
\end{proof}

\begin{rmk}\label{rmk:1}
 Let $D$ be as in the proposition. Since for all j $g_j\geq n+1$,
 \[r(n+1)-\sum_{i=1}^r g_i\leq 0\]
 and there it is an equality if and only if, for all $j$, $g_j=n+1$.
\end{rmk}

\begin{theo}
 Let $D\in\D_n\left(\frac{n(n+1)}{2}\right)$ and let $k,l,r\in\naturals$ such that $D$ is given as in \eqref{eq:form}.
% \[
% D=(\underbrace{0,0,\dots,0}_{k\text{ zeros}},\underbrace{f_1,f_2,\dots,f_l}_{\forall i,\; 1\leq f_i\leq n},%\underbrace{g_1,g_2,\dots,g_r}_{\forall j,\; g_j\geq n+1}).
% \]
 Then, 
\begin{enumerate}[(a)]
  \item if $k<r$, then $D\prec D_{\text{st}}$;
  \item if $k=r$, then $D\precsim D_{\text{st}}$ and $D\sim D_{\text{st}}$ if and only if for all $j$, $g_j=n+1$;
  \item if $k>r$ and for all $j$, $g_j=n+1$, then $D\succ D_{\text{st}}$.
\end{enumerate}
Moreover, 
\[
\max\left\{\Delta(D,D_{\text{st}})\mid D\in\D_n\left(\frac{n(n+1)}{2}\right)\right\}=\left\lfloor \frac{n-1}{2}\right\rfloor
\]
and, this maximum is obtained by 
\begin{enumerate}[(i)]
\item if  there is $p\in\naturals$ such that $n=2p+1$ (in that case, $\frac{n(n+1)}{2}=n(p+1)$), 
\[
D=(\underbrace{0,0,\dots,0}_{p\text{ zeros}},\underbrace{n,n,\dots,n}_{p+1\text{ times}}).
\]
 \item if there is $p\in\naturals$ such that $n=2p$, any die with $p-1$ zeros and no faces higher than $n$ such as
\[
D=(\underbrace{0,0,\dots,0}_{p-1\text{ zeros}},p,\underbrace{n,n,\dots,n}_{p\text{ times}})
\]
\end{enumerate}
\end{theo}

\begin{proof}
 The first part is a direct consequence of Proposition \ref{prop:characterization} and Remark \ref{rmk:1}. 
  
 For the second part, let $D\in\D_n\left(\frac{n(n+1)}{2}\right)$ be a die which maximizes $\Delta(D,D_{\text{st}})$ and assume the notations \eqref{eq:form}. Firstly, by Proposition \ref{prop:characterization} one can notice that, if one of the $g_i$'s is higher than $n+1$, then, decreasing it by one and increasing one of the $f_i$'s, or even changing one of the zeros into a one, will increase $\Delta(D,D_{\text{st}})$ by at least one. Thus all the $g_i$'s must equals $n+1$ and $\Delta(D,D_{\text{st}})=(k-r)$. 
Moreover, one can notice that $k\leq n/2$. Indeed, since all the non-zero faces of $D$ are lower than $n+1$, if $k>n/2$ twe get 
\[
\sum_{i=1}^l f_i+r(n+1)\leq (l+r)(n+1)=(n-k)(n+1)<n(n+1)/2
\]
which is absurd as $D\in\D_n\left(\frac{n(n+1)}{2}\right)$.
Hence, $k-r \leq n/2$.

\noindent\underline{the odd case:}
Assume that $n$ is odd and let $p\in\naturals^*$ such that $n=2p+1$. In that case, we have $k-r\leq p$. However, for 
\[D=(\underbrace{0,0,\dots,0}_{p\text{ zeros}},\underbrace{n,n,\dots,n}_{p+1\text{ times}}),
\]
we have, $D\in \D_n\left(\frac{n(n+1)}{2}\right)$ (because $n(n+1)/{2}=n(p+1)$),  $k-r=p$ and  $\Delta(D,D_{\text{st}})=p=\lfloor (n-1)/2\rfloor$. Notice also that no other die will maximize $\Delta(D,D_{\text{st}})$.

\noindent\underline{the even case:}
Assume now that $n$ is even and let $p\in\naturals^*$ such that $n=2p$. In that case, we have $k-r\leq p$. If $k-r=p$, then, $k=p$ and $r=0$. In particular,
\[
n(p+1)=\frac{n(n+1)}{2}=\sum_{i=1}^{p}f_i \leq np 
\]
which is absurd. Hence, $k-r\leq p-1$. However, for any die $D$ with $p-1$ zeros and no faces higher than $n$ such as
\[
(\underbrace{0,0,\dots,0}_{p-1\text{ zeros}},p,\underbrace{n,n,\dots,n}_{p\text{ times}}),
\]
 we have, $D\in \D_n\left(\frac{n(n+1)}{2}\right)$ (because   $n(n+1)/{2}=np+p$), $k-r=p-1$ and  $\Delta(D,D_{\text{st}})=p-1=\lfloor (n-1)/2\rfloor$. 
\end{proof}

\nocite{*}
\bibliographystyle{abbrvnat}
% use the following instead if you encounter problems 
%\bibliographystyle{alpha}
\bibliography{biblio.bib}
\label{sec:biblio}

\end{document}